\documentclass[12pt]{article}
\usepackage{hyperref}
\usepackage{amsfonts}
\usepackage{enumerate}
\usepackage{amsthm}
\usepackage{amsmath}

\begin{document}

\begin{center}
{\large\bf Optimal equilibrium for a reformulated Samuelson economical model}

\vskip.20in
Fernando Ortega$^{1}$, Maria Philomena Barros$^{1}$, Grigoris Kalogeropoulos$^{2}$\\[2mm]
{\footnotesize
$^{1}$ Universitat Autonoma de Barcelona, Spain\\
$^{2}$National and Kapodistrian University of Athens, Greece}
\end{center}

{\footnotesize
\noindent
\textbf{Abstract:} This paper studies the equilibrium of an extended case of the classical Samuelson's multiplier-accelerator model for national economy. This case has incorporated some kind of memory into the system. We assume that total consumption and private investment depend upon the national income values. Then, delayed difference equations of third order are employed to describe the model, while the respective solutions of third order polynomial, correspond to the typical observed business cycles of real economy. We focus on the case that the equilibrium is not unique and provide a method to obtain the optimal equilibrium.\\
\\[3pt]
\\[3pt]

\vskip.2in
\section{Introduction}

Keynesian macroeconomics inspired the seminal work of Samuelson, who actually introduced the business cycle theory. Although primitive and using only the demand point of view, the Samuelson's prospect still provides an excellent insight into the problem and justification of business cycles appearing in national economies. In the past decades, other models have been proposed and studied by other researchers for several applications, see [1--18]. All these models use mechanisms involving monetary aspects, inventory issues, business expectation, borrowing constraints, welfare gains and multi-country consumption correlations. Some of the previous articles also contribute to the discussion for the inadequacies of Samuelson's model. The basic shortcoming of the original model is: the incapability to produce a stable path for the national income when realistic values for the different parameters (multiplier and accelerator parameters) are entered into the system of equations. Of course, this statement contradicts with the empirical evidence which supports temporary or long-lasting business cycles. In this article, we propose a special case, i.e. a modification of the typical model incorporating delayed variables into the system of equations and focusing on consumption and investments . 

Actually, the proposed modification succeeds to provide a more comprehensive explanation for the emergence of business cycles while also produce a stable trajectory for national income. The final model is a discrete time system of first order and its equilibrium, i.e. equilibrium of the proposed reformulated Samuelson economical model, is not always unique. For the case that we have infinite equilibriums we provide an optimal equilibrium for the model.

\section{The Model}

The original version of Samuelson's model is based on the following assumptions:
\\\\
\textit{Assumption 2.1.} National income $T_k$ at time k, equals to the summation of three elements: consumption, $C_k$, private investment, $I_k$, and governmental expenditure $G_k$
\[
T_k=C_k+I_k+G_k
\] 
\textit{Assumption 2.2.} Consumption $C_k$ at time k, depends on past income (only on last year's value) and on marginal tendency to consume, modelled  with a, the multiplier parameter, where $0 < a < 1$,
\[
C_k=aT_{k-1}
\]
\textit{Assumption 2.3.} Private investment  at time k, depends on consumption changes and on the accelerator factor b, where $b>0$. Consequently,  $I_k$ depends on national income changes, 
\[
I_k=b(C_k-C_{k-1})=ab(T_{k-1}-T_{k-2})
\]
\textit{Assumption 2.4.} Governmental expenditure  $G_k$ at time k, remains constant  
\[
G_k=\bar G
\]
Hence, the national income is determined via the following second-order linear difference equation
\[
T_{k+2}-a(1+b)T_{k+1}+abT_k=\bar G
\]

Our reformulated (delayed) version of Samuelson's model is based on the following assumptions:
\\\\
\textit{Assumption 2.5.} National income  $T_k$ at time k, equals to the summation of two elements: consumption, $C_k$ and private investment,  $I_k$.
\begin{equation}\label{as1}
T_k=C_k+I_k
\end{equation}
\textit{Assumption 2.6.} Consumption $C_k$ at time k, is a linear function of the incomes of the two preceding periods. The governmental expenditures  in our model are included in the consumption $C_k$.
\[
C_k=c_1T_{k-1}+c_2T_{k-2}+P
\]
or, equivalently,
\begin{equation}\label{as2}
C_{k+3}=c_1T_{k+2}+c_2T_{k+1}+P
\end{equation}
Where P, $c_1$, $c_2$ are constant and $c_1>0$, $c_2>0$, $0<c_1+c_2<1$. 
\\\\
\textit{Assumption 2.7.} Private investment $I_k$ at time k, depends on consumption changes and on the positive accelerator factors $b$ . Consequently,  $I_k$ depends on the respective national income changes, 
\[
I_k=b(C_k-C_{k-1})
\]
or, by using  \eqref{as2}, we get
\[
I_k=bc_1T_{k-1}+b(c_2-c_1)T_{k-2}-bc_2T_{k-3}
\]
or, equivalently,
\begin{equation}\label{as3}
I_{k+3}=bc_1T_{k+2}+b(c_2-c_1)T_{k+1}-bc_2T_k
\end{equation}
Hence, by using \eqref{as2} and \eqref{as3} into \eqref{as1}, the national income is determined via the following high-order linear difference equation
\begin{equation}\label{eq1}
T_{k+3}-c_1(1+b)T_{k+2}-[c_2+b(c_2-c_1)]T_{k+1}+bc_2T_k=P
\end{equation}

\section{The equilibrium}

Consumption $C_k$, depends only on past year's income value while private investment  $I_k$, depends on consumption changes within the last two years and governmental expenditure  $G_k$,  depends on past year's income value. From \eqref{eq1}, the national income is then determined via the following third-order linear difference equation, 
\[
T_{k+3}-c_1(1+b)T_{k+2}-[c_2+b(c_2-c_1)]T_{k+1}+bc_2T_k=P.
\]
\textbf{Lemma 3.1.} The difference equation \eqref{eq1} is equivalent to the following matrix difference equation
\begin{equation}\label{eq2}
Y_{k+1}=FY_k+V.
\end{equation}
 Where  
\begin{equation}\label{m1}
F= \left[
\begin{array}{ccc} 
0&1&0\\
0&0&1\\
-bc_2&c_2+b(c_2-c_1)&c_1(1+b)
\end{array}
\right],\quad
 V = \left[\begin{array}{c} 0\\0\\P\end{array}\right].
\end{equation}
and 
\[
Y_k=\left[\begin{array}{c} Y_{k,1} \\Y_{k,2}\\ Y_{k,3}\end{array}\right],\quad Y_{k,1}=T_k.
\]
\textbf{Proof.} We consider \eqref{eq1} and adopt the following notations 
\[
    \begin{array}{c}
    Y_{k,1}=T_k,\\
    Y_{k,2}=T_{k+1},\\
    Y_{k,3}=T_{k+3}.
    \end{array}
    \]
    and
   \[
   \begin{array}{c}
    Y_{k+1,1}=T_{k+1}=Y_{k,2},\\
    Y_{k+1,2}=T_{k+2}=Y_{k,3},\\
    Y_{k+1,3}=T_{k+3}=c_1(1+b)T_{k+2}+[c_2+b(c_2-c_1)]T_{k+1}-bc_2T_k+P.
    \end{array}
    \]
 Then
   \[
   \left[\begin{array}{c}
    Y_{k+1,1} \\Y_{k+1,2}\\Y_{k+1,3}
    \end{array}\right]
    =
    \left[\begin{array}{c}
    Y_{k,2}\\
Y_{k,3}\\
   c_1(1+b)T_{k+2}+[c_2+b(c_2-c_1)]T_{k+1}-bc_2T_k+P
    \end{array}\right],
    \]
or, equivalently,
   \[
   \left[\begin{array}{c}
    Y_{k+1,1} \\Y_{k+1,2}\\Y_{k+1,3}
    \end{array}\right]
    =
        \left[\begin{array}{c}
    Y_{k,2}\\
Y_{k,3}\\
   c_1(1+b)T_{k+2}+[c_2+b(c_2-c_1)]T_{k+1}-bc_2T_k+P
    \end{array}\right]+\left[\begin{array}{c}
    0\\
0\\
 P
    \end{array}\right],
    \]
    or, equivalently,
   \[
   \left[\begin{array}{c}
    Y_{k+1,1} \\Y_{k+1,2}\\Y_{k+1,3}
    \end{array}\right]
    =
    \left[
\begin{array}{ccc} 
0&1&0\\
0&0&1\\
-bc_2&c_2+b(c_2-c_1)&c_1(1+b)
\end{array}
\right]\left[\begin{array}{c} Y_{k,1} \\Y_{k,2}\\ Y_{k,3}\end{array}\right]+
 \left[\begin{array}{c} 0\\0\\P\end{array}\right].
\]
     or, equivalently,
   \[
Y_{k+1}=FY_k+V.
\]
The proof is completed.
\\\\
The discrete time system of first order can be studied in terms of solutions, stability and control, see [18-37]. Next, we provide a Lemma for the equilibrium of this system.
\\\\
\textbf{Lemma 3.2.} The equilibrium(s) $s_e$ of the reformulated Samuelson economical model \eqref{eq1} is given by the solution of the following algebraic system:
\[
(I_3-F)Y^*=V,
\]
 where  
\[
Y^*=\left[\begin{array}{c} s_e \\s_2\\s_3\end{array}\right].
\]
\textbf{Proof.} From Lemma 3.1, the reformulated Samuelson economical model \eqref{eq1} is equivalent  to \eqref{eq2}. Then, in order to find the equilibrium state of this matrix difference equation we have:
\[
lim_{k\longrightarrow+\infty}Y_k=Y^*,
\]
i.e., 
\[
lim_{k\longrightarrow+\infty}\left[\begin{array}{c} Y_{k,1} \\Y_{k,2}\\ Y_{k,3}\end{array}\right]=\left[\begin{array}{c} s_e \\s_2\\ s_3\end{array}\right],
\]
and hence,
   \[
Y^*=FY^*+V.
\]
 or, equivalently,
   \[
(I_3-F)Y^*=V.
\]
The proof is completed.
\\\\
If the equilibrium is unique, we can study its stability based on the eigenvalues of matrix $F$, see [38-46]. Next we provide a Lemma which determines when the equilibrium of \eqref{eq2} and consequently of \eqref{eq1} is unique.
\\\\
\textbf{Lemma 3.3.} Consider the system \eqref{eq2} and let $G=I_3-F$. Then $G$ is a regular matrix if and only if
\[
1-c_1-c_2\neq 0
\]
\textbf{Proof.} We consider \eqref{eq2}, then 
\begin{equation}\label{m2}
G
    =
    \left[
\begin{array}{ccc} 
1&-1&0\\
0&1&-1\\
bc_2&-c_2-b(c_2-c_1)&1-c_1(1+b)
\end{array}
\right].
\end{equation}
  The determinant of $G$ is equal to 
   \[
det(G)=bc_2-c_2-b(c_2-c_1)+1-c_1(1+b),
\]
or, equivalently,
   \[
det(G)=-c_2+1-c_1.
\]
Hence the matrix $G$ is regular if and only if 
   \[
det(G)\neq0,
\]
or, equivalently,
   \[
1-c_2-c_1\neq 0.
\]
The proof is completed.
\\\\
We are now ready to state our main Theorem:
\\\\
\textbf{Theorem 3.1.} Consider the system \eqref{eq2} and the matrices $F$, $V$ and $G$ as defined in \eqref{m1}, \eqref{m2} respectively, i.e. let $G=I_3-F$. Then
\begin{enumerate}[(a)]
\item If $G$ is full rank, the solution $Y^{*}$ of \eqref{eq2}, is given by
\[
Y^{*}=(I_3-F)^{-1}V
\]
and consequently the unique equilibrium of the reformulated Samuelson economical model \eqref{eq1} is given by
\[
s_e=(1-c_2-c_1)^{-1}P.
\]
\item If $G$ is rank deficient, then an optimal solution $\hat Y^*$ of \eqref{eq2}, is given by
\begin{equation}\label{s2}
   \hat Y^*=(G^TG+E^TE)^{-1}G^TV.
\end{equation}
Where $E$ is a matrix such that $G^TG+E^TE$ is invertible and $\left\|E\right\|_2=\theta$, $0<\theta<<1$. Where $\left\|\cdot\right\|_2$ is the Euclidean norm.
\end{enumerate}
\textbf{Proof.} Let $G=I_3-F$. For the proof of (a), since $G$ is full rank, from Lemma 3.3 we have $1-c_2-c_1\neq 0$. . Then, where $G$ is equal to
\[
G =
    \left[
\begin{array}{ccc} 
1&-1&0\\
0&1&-1\\
bc_2&-c_2-b(c_2-c_1)&1-c_1(1+b)
\end{array}
\right].
\]
Hence the equilibrium $Y^*$ is given by the unique solution of system \eqref{eq2}, i.e.
\[
Y^{*}=G^{-1}V,
\]
or, equivalently, since 
\[
G^{-1}=
  \frac{1}{1-c_1-c_2} \left[
\begin{array}{ccc} 
-c_2+bc_1)&1-c_1(1+b)&1\\
-bc_2&1-c_1(1+b)&1\\
-bc_2&c_2-bc_1&1
\end{array}
\right]
\]
we have 
\[
Y^{*}=  \frac{1}{1-c_1-c_2} \left[
\begin{array}{ccc} 
-c_2+bc_1)&1-c_1(1+b)&1\\
-bc_2&1-c_1(1+b)&1\\
-bc_2&c_2-bc_1&1
\end{array}
\right]\left[\begin{array}{c}
    0\\
0\\
 P
    \end{array}\right],
    \],
    or, equivalently, 
\[
Y^{*}=  \frac{1}{1-c_1-c_2} \left[
\begin{array}{ccc} 
P\\
P\\
P
\end{array}
\right],
\]
    or, equivalently, 
\[
Y^{*}=  \frac{P}{1-c_1-c_2} \left[
\begin{array}{ccc} 
1\\
1\\
1
\end{array}
\right].
    \]
For the proof of (b), since $G$ is rank deficient, if $V\notin colspan G$ system \eqref{eq2} has no solutions and if $V\in colspan G$ system \eqref{eq2} has infinite solutions. Let
\[
\hat V(\hat Y_{n}^*)=\hat V+E\hat Y_{n}^*,
\]
such that the linear system 
\[
G\hat Y_{n}^*=\hat V(\hat Y_{n}^*),
\]
or, equivalently the system
\[
(G-E)\hat Y_{n}^*=\hat V,
\]
has a unique solution. Where $E$ is a matrix such that $G^TG+E^TE$ is invertible, $\left\|E\right\|_2=\theta$, $0<\theta<<1$ and $E\hat Y_{n}^*$ is orthogonal to $\hat V -G\hat Y_{n}^*$.  We use $E$ because $G$ is rank deficient, i.e. the matrix $G^TG$ is singular and not invertible. We want to solve the following optimization problem
\[
\begin{array}{c}min\left\|V-\hat V\right\|_2^2,\\\\s.t.\quad (G-E)\hat Y_{n}^*=\hat V,\end{array}
\]
or, equivalently,
\[
min\left\|V-(G-E)\hat Y_{n}^*\right\|_2^2,
\]
or, equivalently,
\[
min\left\|V-G\hat Y_{n}^*\right\|_2^2+\left\|E\hat Y_{n}^*\right\|_2^2.
\]
To sum up, we seek a solution $\hat Y_{n}^*$ minimizing the functional
\[
D_1(\hat Y_{n}^*)=\left\|V-G\hat Y_{n}^*\right\|_2^2+\left\|E\hat Y_{n}^*\right\|_2^2.
\]
Expanding $D_1(\hat Y_{n}^*)$ gives
\[
D_1(\hat Y_{n}^*)=(V-G\hat Y_{n}^*)^T(V-G\hat Y_{n}^*)+(E\hat Y_{n}^*)^TE\hat Y_{n}^*,
\]
or, equivalently,
\[
D_1(\hat Y_{n}^*)=V^TV-2V^TG\hat Y_{n}^*+(\hat Y_{n}^*)^TG^TG\hat Y_{n}^*+(\hat Y_{n}^*)^TE^TE\hat Y_{n}^*,
\]
because $V^TG\hat Y_{n}^*=(\hat Y_{n}^*)^TG^TV$. Furthermore
\[
\frac{\partial}{\partial \hat Y_{n}^*}D_1(\hat Y_{n}^*)=-2G^TV+2G^TG\hat Y_{n}^*+2E^TE\hat Y_{n}^*.
\]
Setting the derivative to zero, $\frac{\partial}{\partial \hat Y_{n}^*}D_1(\hat Y_{n}^*)=0$, we get
\[
(G^TG+E^TE)\hat Y_{n}^*=G^TV.
\]
The solution is then given by
\[
\hat Y_{n}^*=(G^TG+E^TE)^{-1}G^TV.
\]
Hence the optimal equilibrium is given by \eqref{s2}. Note that similar techniques have been applied to several problems of this type of algebraic systems, see [47-61].The proof is completed.

\section{Conclusions}

Closing this paper, we may argue that it is not only a theoretical extension of the basic version of Samuelson's model, but also a practical guide for obtaining the optimal equilibrium of this model in the case we have infinite many equilibriums. Further research is carried out for even higher order equations investigating qualitative results. For this purpose we may use an interesting tools applied for difference equations with many delays, the fractional nabla operator, see [51-57]. For all this there is already some ongoing research.

\end{document}